\documentclass{siamltex}
\usepackage{amsfonts,latexsym,amssymb}
\usepackage{graphicx,amsmath}
\usepackage{algorithm}	
\usepackage{algorithmic}	
\newcommand{\R}{{\mathbb R}}
\newcommand{\N}{{\mathbb N}}

\newcommand{\bx}{\mbox{\boldmath{$x$}}}

\newcommand{\bw}{\mbox{\boldmath{$w$}}}
\newcommand{\be}{\mbox{\boldmath{$e$}}}

\newcommand{\bg}{\mbox{\boldmath{$g$}}}

\newcommand{\by}{\mbox{\boldmath{$y$}}}
\newcommand{\bz}{\mbox{\boldmath{$z$}}}
\newcommand{\sby}{\mbox{\boldmath{${\scriptstyle y}$}}}
\newcommand{\sbx}{\mbox{\boldmath{${\scriptstyle x}$}}}

\newcommand{\sbz}{\mbox{\boldmath{${\scriptstyle z}$}}}

\newcommand{\bzero}{\mbox{\boldmath{$0$}}}
\newcommand{\bl}{\begin{list}{ \ }{
\leftmargin=.325in}}
\newcommand{\el}{\end{list}}

\newcommand{\ri}{\mathbf{i}}
\newcommand{\ci}{\mathbf{j}}

\newcommand{\norm}[1]{\mbox{$\parallel\!#1\!\parallel$}}

\usepackage{booktabs}
\usepackage{xcolor}%
\definecolor{SPECgreen}{rgb}{0,0.75,0}%
\definecolor{SPECorange}{rgb}{1.0,.5625,0}%
\definecolor{SPECred}{rgb}{0.75,0,0}

\definecolor{SPECblue}{rgb}{0,0,0.75}
\definecolor{SPECblack}{rgb}{0.25,0.25,0.25}
\definecolor{SPECgray}{rgb}{0.625,0.625,0.625}
\definecolor{citecol}{rgb}{0.75,0,0}
\definecolor{linkcol}{rgb}{0,0,0.75}
\definecolor{SPECpink}{rgb}{1.0,0.1,0.6}
\definecolor{SPECorangeg}{rgb}{0.609375,.609375,0.609375}
\definecolor{SPECblueg}{rgb}{0.0525,0.0525,0.0525}

\begin{document}

\date{}
\title{Adaptive cross approximation for Tikhonov regularization in general form}

\author{
T. Mach\thanks{Institute for Mathematics, University of Potsdam, 14476 Potsdam, Germany. 
E-mail: \texttt{mach@uni-potsdam.de}.} 
\and
L. Reichel\thanks{Department of Mathematical Sciences, Kent State University, Kent, OH
44242, USA. E-mail: \texttt{reichel@math.kent.edu}.}
\and
M. Van Barel\thanks{Department of Computer Science, KU~Leuven, Celestijnenlaan 200A, 3001 
Leuven (Heverlee), Belgium. E-mail: \texttt{marc.vanbarel@cs.kuleuven.be}.}
}
\maketitle

\begin{center}
\small{Dedicated to Claude Brezinski on the occasion of his 80th birthday.}
\end{center}

\begin{abstract}
Many problems in Science and Engineering give rise to linear integral equations of the 
first kind with a smooth kernel. Discretization of the integral operator yields a 
matrix, whose singular values cluster at the origin. We describe the approximation of 
such matrices by adaptive cross approximation, which avoids forming the entire matrix. The
choice of the number of steps of adaptive cross approximation is discussed. The 
discretized right-hand side represents data that commonly are contaminated by measurement 
error. Solution of the linear system of equations so obtained is not meaningful because 
the matrix determined by adaptive cross approximation is rank-deficient. We remedy this 
difficulty by using Tikhonov regularization and discuss how a fairly general 
regularization matrix can be used. Computed examples illustrate that the use of a 
regularization matrix different from the identity can improve the quality of the computed 
approximate solutions significantly.
\end{abstract}


\section{Introduction}\label{sec1}
Linear integral equations of the first kind,
\begin{equation}\label{inteq}
\int_{\Omega_1} \kappa(s,t) x(t)dt = g(s),\qquad s\in\Omega_2,
\end{equation}
with a smooth kernel $\kappa$ arise in many applications, including remote sensing, 
computerized tomography, and image restoration. Here $\Omega_i$ denotes a subset 
of $\R^{d_i}$ for some positive integer $d_i$. The solution of \eqref{inteq} is an 
\emph{ill-posed problem}, because the singular values of the integral operator cluster at 
the origin; see, e.g., \cite{EHN,Gr} for introductions to ill-posed problems. 

Discretization of \eqref{inteq} yields a linear system of equations
\begin{equation}\label{linsys}
A\bx=\bg, \qquad A \in {\R}^{n\times n}, \qquad \bg \in {\R}^n,
\end{equation}
with a matrix, whose singular values coalesce at the origin. This makes the matrix $A$
severely ill-conditioned and possibly rank-deficient; we measure the conditioning of a 
matrix with its condition number, which is the ratio of the largest and smallest singular 
values. Linear systems of equations with a matrix of this kind are often referred to as 
\emph{discrete ill-posed problems}; see, e.g., \cite{Ha1}. We will for notational 
simplicity assume the matrix $A$ to be square, however, the method described also can be 
applied, after minor modifications, when $A$ is rectangular, in which case the linear 
system of equations \eqref{linsys} is replaced by a least-squares problem.

In many applications, the right-hand side vector $\bg$ in \eqref{linsys} represents 
measured data and is contaminated by a measurement error $\be$. Due to the severe 
ill-conditioning of $A$, straightforward solution of \eqref{linsys} typically yields a 
computed solution that is severely contaminated by propagated error, and therefore is not
useful. To circumvent this difficulty, the linear system of equations \eqref{linsys} 
commonly is replaced by a nearby problem, whose solution is less sensitive to the error 
$\be$ in $\bg$. This replacement is referred to as \emph{regularization}. Tikhonov 
regularization is possibly the most popular and well understood regularization method. It
replaces the linear system of equations \eqref{linsys} by a penalized least-squares 
problem of the form
\begin{equation}\label{tikhonov}
\min_{\sbx\in\R^n}\{\|A\bx-\bg\|^2+\mu\|L\bx\|^2\},
\end{equation}
where $L\in\R^{p\times n}$ is referred to as the regularization matrix and $\mu>0$ as the 
regularization parameter. The problem \eqref{tikhonov} is said to be in 
\emph{standard form} when $L$ is the identity; otherwise \eqref{tikhonov} is in 
\emph{general form}. Throughout this paper $\|\cdot\|$ denotes the Euclidean vector 
norm or the spectral matrix norm. 

The choice of regularization parameter $\mu$ is important for the quality of the 
computed solution: a too small value results in a computed solution that is contaminated 
by needlessly much propagated error, while a too large value yields an unnecessarily 
smooth solution that may lack details of interest. Generally, a suitable value of $\mu$
is not known a priori, but has to be computed during the solution process. This typically
requires that \eqref{tikhonov} be solved for several $\mu$-values. Methods for determining
a suitable value of $\mu$ include the L-curve criterion, generalized cross validation, and
the discrepancy principle; see, e.g., \cite{BRS,BRZRS,FRRS,KR,PRRY,RR} for discussions of 
properties of these and other methods. 

The matrix $L$ is assumed to be chosen so that
\begin{equation}\label{nullcond}
{\mathcal N}(A)\cap{\mathcal N}(L)=\{\bzero\},
\end{equation}
where ${\mathcal N}(M)$ denotes the null space of the matrix $M$. Then the Tikhonov 
minimization problem \eqref{tikhonov} has the unique solution
\begin{equation}\label{tiksol}
\bx_\mu:=(A^TA+\mu L^TL)^{-1}A^T\bg
\end{equation}
for any $\mu>0$; see, e.g, \cite{Ha1} for details. Here and below the superscript $^T$ 
denotes transposition. We are interested in the situation when the matrices $A$ and $L$
are so large that it is impossible or undesirable to compute the solution \eqref{tiksol}
by Cholesky factorization of the matrix $A^TA+\mu L^TL$. In fact, we would like to 
avoid evaluating all the entries of $A$. We describe how $A$ can be approximated by a much
smaller matrix, without evaluating all matrix entries, by applying adaptive cross 
approximation. 

The application of cross approximation to matrices that stem from the discretization of 
Fredholm integral equations of the second kind has received considerable attention in the
literature, see, e.g., \cite{q470,q699,FVB,GTZ,Ty}; however, the use of cross 
approximation in the context of solving linear discrete ill-posed problems has not been
thoroughly studied. 

The use of adaptive approximation to the approximate solution of \eqref{tikhonov} when $L$
is the identity is discussed in \cite{MRVV}. This paper extends this discussion to general
regularization matrices $L$.
Our interest in this extension of the method in \cite{MRVV} stems from the fact
that the use of a suitably chosen regularization matrix $L$ can deliver solutions of 
higher quality than $L=I$; see, e.g., \cite{HNR,NR,RYe} for illustrations and discussions 
on various ways of constructing regularization matrices. Roughly, $L$ should be chosen so 
as not to damp known important features of the desired solution, while damping the 
propagated error stemming from the error in $\bg$. 

In the computed examples of Section \ref{sec4b}, we use the \emph{discrepancy principle} 
to determine $\mu>0$. Let $\widehat{\bg}\in\R^n$ denote the unknown error-free vector 
associated with the right-hand side $\bg$ in \eqref{linsys}, i.e., 
$\bg=\widehat{\bg}+\be$. Assume that the linear system of equations with the error-free 
right-hand side, 
\begin{equation}\label{exsys}
A\bx=\widehat{\bg},
\end{equation}
is consistent and that a fairly accurate bound $\|\be\|\leq\delta$ is known. The 
discrepancy principle prescribes that the regularization parameter $\mu>0$ be determined 
so that the Tikhonov solution \eqref{tiksol} satisfies
\begin{equation}\label{discrp}
\|A\bx_\mu-\bg\|=\eta\delta, 
\end{equation}
where $\eta>1$ is a user-specified parameter that is independent of $\delta$. It can be
shown that when $\delta$ tends to zero, $\bx_\mu$ converges to the minimal-norm solution,
$\widehat{\bx}$, of \eqref{exsys}; see, e.g., \cite{EHN} for a proof in a Hilbert space 
setting. We remark that the determination of $\mu>0$ such that $\bx_\mu$ satisfies
\eqref{discrp} typically requires the solution of \eqref{tikhonov} for several 
$\mu$-values.

The present paper is concerned with the situation when the matrix $A\in\R^{n\times n}$ in
\eqref{linsys} is large. Then the evaluation of all entries of $A$ can be quite 
time-consuming. Cross approximation, also known as skeleton approximation, of $A$ reduces 
this time by approximating $A$ by a matrix $M_k\in\R^{n\times n}$ that consists of 
$k\ll n$ rows and columns of $A$. We would like to choose the rows and columns of $A$ so 
that $M_k$ approximates $A$ well and is easy to compute with. 

This paper is organized as follows. Section \ref{sec2} reviews the application of adaptive
cross approximation to the approximation of $A$ by a matrix of low rank. In Section 
\ref{sec3}, we describe the application of adaptive cross approximation to the 
approximation of the Tikhonov equation \eqref{tikhonov}. Section \ref{sec4b} reports a few
computed examples, and concluding remarks can be found in Section \ref{sec5}.

We conclude this section with an example that leads to a large linear discrete ill-posed
problem, whose solution is difficult to compute using straightforward discretization.

Example 1.1. Consider the solution of the Fredholm integral equation of the first
kind 
\begin{equation}\label{potential}
\int_{S}\frac{\sigma(\by)}{4\pi\epsilon_{0}\norm{\bx-\by}}\,d\by=\phi(\bx),\qquad\bx\in S,
\end{equation}
where $\phi$ is a given electric potential, $S$ is a surface with electrodes, 
$\sigma(\by)$ denotes the density of the charge on $S$, and $\epsilon_{0}$ stands for the
electric permittivity in vacuum. We would like to determine $\sigma$ from $\phi$, and 
assume that $\phi$ is chosen so that \eqref{potential} has a solution. The computation 
of a solution of \eqref{potential} is an ill-posed problem. Using a weak formulation and 
discretization lead to a dense symmetric matrix $K=[k_{ij}]_{i,j=1}^n\in\R^{n\times n}$ 
with entries
\[
k_{ij} = \int_{S}\int_{S} \frac{v_{j}(\bx)v_{i}(\by)}
  {4\pi\epsilon_{0}\norm{\bx-\by}}\,d\by\,d\bx,\quad i,j=1,2,\ldots,n.
\]
This matrix can be expensive to store and handle when the discretization is fine.
Employing a hierarchical compression with $\mathcal{H}^{2}$-matrices reduces the required 
storage to $O(n)$, with a large constant hidden in the $O(\cdot)$, and allows 
matrix-vector product evaluations in $O(n)$ arithmetic floating point operations (flops);
see \cite{BoermH2}.

A fine discretization with $n=262,146$ nodes results in a large, $n\times n$, dense
matrix. We used the H2Lib library \cite{H2lib} for the computations and base this example 
on one of the standard examples provided in this library. Without compression, 512 GB of memory are 
needed to store the matrix. On a laptop computer with an Intel Core i710710U
CPU and 16 GB of RAM it took 1103 s to assemble the matrix $K$ in the compressed
$\mathcal{H}^{2}$-matrix format. The matrix required 15.45 GB of storage, thus almost
all the available RAM. Carrying out one matrix-vector product evaluation required 1596 s,
that is 44\% more time than for assembling the matrix. The reason for this is that the
$O(n)$ flops require a significant amount of communication between faster and slower
storage. This example illustrates that there are linear discrete ill-posed problems of 
interest that are difficult to solve on a laptop computer, even if a significant amount of
memory is available. It therefore is important to develop methods that are able to 
determine approximations of dense matrices that requires less computer storage and 
less CPU time for the evaluation of matrix-vector products. \hfill$\Box$

\section{Adaptive cross approximation}\label{sec2}
Cross approximation of a large matrix $A\in\R^{n\times n}$ determines an approximation 
$M_k\in\R^{n\times n}$ of rank at most $k\ll n$. All entries of $M_k$ can be evaluated 
much faster than all entries of $A$, because $M_k$ is constructed from only $k$ rows and
columns of $A$. We would like $M_k$ to be an accurate approximation of $A$. This is 
achieved by a careful choice of the $k$ rows and columns of $A$ that define the matrix
$M_k$. A cross approximation method is said to be \emph{adaptive} when the rows and 
columns of $A$ that determine $M_k$ (and $k$) are chosen depending on properties of $A$ 
revealed during the computations; see \cite{q467,q699}. 

We outline the adaptive cross approximation method for a general square nonsymmetric 
matrix described in \cite{MRVV}. This method is an adaptation of the scheme in \cite{FVB}
to the approximation of the matrix of linear discrete ill-posed problems. When $A$ is 
symmetric, the matrix $M_k$ can be chosen to be symmetric. This roughly halves the storage 
requirement for $M_k$. Both the situations when $A$ is symmetric positive definite or 
symmetric indefinite are discussed in \cite{MRVV}. We therefore will not dwell on these
special cases in the present paper.

Let the matrix $A\in\R^{n\times n}$ be nonsymmetric and choose $k$ rows of $A$ with 
indices in $\N=\{1,2,\dots,n\}$. We let $\ri,\ci\in\N^k$ denote index vectors with $k$
entries in $\N$. The submatrices $A_{(\ri,:)}$ and $A_{(:,\ci)}$ of $A$ are made up of the 
$k$ rows with indices $\ri$ and the $k$ columns with indices $\ci$, respectively. 
Moreover, the \emph{core matrix} $A_{(\ri,\ci)}$ is made up of $k$ rows and columns of $A$.
Assume that this matrix is nonsingular. Then the rows and columns of the matrix
\[
  M_k = A_{(:,\ci)} A_{(\ri,\ci)}^{-1} A_{(\ri,:)}
\]
are equal to the corresponding rows and columns of $A$; when $A$ is of rank $k$, we have 
$M_k=A$. 

Goreinov et al. \cite{q468} show that it is beneficial to choose the index vectors $\ri$ 
and $\ci$ so that $A_{(\ri,\ci)}$ is a submatrix of $A$ of maximal volume, i.e., so that
the modulus of the determinant of $A_{(\ri,\ci)}$ is maximal. However, it is difficult to
determine such index vectors $\ri$ and $\ci$. We therefore seek to determine a low-rank 
matrix $M_k$ that is a sufficiently accurate approximation of $A$ by a greedy algorithm. 
Suppose that we already have computed an approximation 
\[
M_{k-1} = \sum_{\ell=1}^{k-1}\bw^{(c)}_{\ell}(\bw^{(r)}_\ell)^{T},\qquad
\bw^{(c)}_{\ell},\bw^{(r)}_\ell\in\R^n, 
\]
of rank at most $k-1$ of $A$. To compute the next approximation, $M_k$, of $A$ of rank at
most $k$, we determine a row index $i^*$ and a column index $j^*$ by looking for the index 
of the maximum element in magnitude in the previously computed vectors $\bw^{(c)}_{k-1}$ 
(for index $i^{*}$) and $\bw^{(r)}_{k-1}$ (for index $j^{*}$). The vector 
$(\bw^{(r)}_1)^T$ can be chosen as an arbitrary row of $A$. We will let $(\bw^{(r)}_1)^T$ 
be the first row of $A$ in the computed examples of Section \ref{sec4b}. The vector 
$\bw^{(c)}_{1}$ can be chosen in a similar way.

In the simplest form of cross approximation, the determination of the vectors 
$\bw^{(c)}_k$ and $\bw^{(r)}_k$ only requires the entries in row $i^*$ and column 
$j^*$ of $A$ and the elements of already computed vectors $\bw^{(c)}_{\ell}$ and 
$(\bw^{(r)}_\ell)^{T}$, $\ell=1,2,\ldots,k-1$:
\begin{eqnarray*}
  (\bw^{(r)}_k)_j & = & A_{i_{k}^{*},j}-\sum_{\ell=1}^{k-1}(\bw^{(c)}_\ell)_{i_{k}^{*}}
  (\bw^{(r)}_\ell)_{j}, \hspace{0.5cm} \delta = (\bw^{(r)}_k)_{j^{*}_{k}},\\
  (\bw^{(c)}_k)_{i}& =&\displaystyle\frac{1}{\delta}\Big( A_{i,j_{k}^{*}}-
  \sum_{\ell=1}^{k-1}(\bw^{(c)}_\ell)_{i}(\bw^{(r)}_{\ell})_{j_{k}^{*}}\Big).
\end{eqnarray*}
A new skeleton is obtained from the remainder,
\[
R_{k}=A-\sum_{\ell=1}^{k}\bw^{(c)}_{\ell}(\bw^{(r)}_{\ell})^{T},
\]
without explicitly computing all entries of the matrix $R_k$.

The required number of rank-one matrices, $k$, that make up $M_k$ is generally not known a
priori. We would like the difference $A-M_k$ to be of small norm. However, we cannot
evaluate this difference, because most entries of $A$ are not known. Following \cite{FVB}, 
we include $t$ randomly chosen matrix entries $A_{i_\ell,j_\ell}$, for $\ell=1,2,\dots,t$, 
with $i_\ell,j_\ell\in\N$. Define for future reference the set
\begin{eqnarray}\label{set_pi}
\Pi = \{ (i_\ell,j_\ell),~ \forall~\ell=1,2,\dots,t\}
\end{eqnarray}
When a new skeleton is determined, the values of these entries are updated by Subtraction 
from the available skeletons,
\begin{equation}\label{Rk}
 (R_{k})_{i_\ell,j_\ell}=(R_{k-1})_{i_\ell,j_\ell}-
 (\bw^{(c)}_k)_{i_\ell}(\bw^{(r)}_k)_{j_\ell},
\end{equation}
with $(R_0)_{i_\ell,j_\ell}=A_{i_\ell,j_\ell}$. 
The values $(R_{k})_{ i_{\ell},j_{\ell}}$ are used in subsection \ref{sec3_4} as part of the stopping criterion
to determine the final value $k^*$ for $k$.
The value of $t$ is a percentage of the total number of entries. The choice of 
$t$ should depend on properties of the matrix $A$; see \cite{FVB,MRVV} for further 
details. An algorithm is presented in \cite{MRVV}.

\section{Tikhonov regularization in general form}\label{sec3}
This section discusses how to combine adaptive cross approximation with Tikhonov
regularization in general form.

\subsection{Using adaptive cross approximation}\label{ssec3_1}
The matrix $M_k$, whose computation was outlined in the previous section, is of the form
\begin{equation}\label{Mk}
M_k=W_k^{(c)}(W_k^{(r)})^T,
\end{equation}
where
\[
W_k^{(c)}=[\bw_1^{(c)},\bw_2^{(c)},\ldots,\bw_k^{(c)}]\in\R^{n\times k},\qquad
W_k^{(r)}=[\bw_1^{(r)},\bw_2^{(r)},\ldots,\bw_k^{(r)}]\in\R^{n\times k}.
\]
Compute the skinny QR factorizations
\begin{equation}\label{QRfact}
W_k^{(c)}=Q_k^{(c)}R_k^{(c)},\qquad W_k^{(r)}=Q_k^{(r)}R_k^{(r)},
\end{equation}
where the matrices $Q_k^{(c)},Q_k^{(r)}\in\R^{n\times k}$ have orthonormal columns and the 
matrices $R_k^{(c)},R_k^{(r)}\in\R^{k\times k}$ are upper triangular. The factorizations
\eqref{QRfact} can be computed by the Householder-QR method or by factorization methods 
that are designed to perform efficiently on modern computers, such as the methods 
described in \cite{CPR,Er,YNYF}. 

Combining \eqref{Mk} and \eqref{QRfact} yields 
\begin{equation}\label{QRMk} 
M_k=Q_k^{(c)}R_k^{(c)}(R_k^{(r)})^T (Q_k^{(r)})^T.
\end{equation}
Replacing $A$ by $M_k$ in \eqref{tikhonov} gives the minimization problem 
\begin{equation}\label{tikhonov2}
\min_{\sbx\in\R^n}\{\|M_k\bx-\bg\|^2+\mu\|L\bx\|^2\},
\end{equation}
which can be solved in several ways. If the matrix $L$ has a special
structure, such as being banded with small bandwidth, then it may be attractive to
transform \eqref{tikhonov2} to standard form by a technique described by Eld\'en 
\cite{El}. Regularization matrices $L$ with a small bandwidth arise, e.g., when $L$ 
represents a finite difference approximation of a differential operator in one 
space-dimension. It also is easy to transform \eqref{tikhonov2} to standard form when $L$
is an orthogonal projector; see \cite{MRS}. 

In the remainder of this section, we discuss the situation when $L$ is such that 
transformation of \eqref{tikhonov2} to standard form as described in \cite{El} is too 
expensive to be attractive. This is the case, for instance, when $L$ represents a finite 
difference approximation of a differential operator in two or more space-dimensions. This 
kind of matrices $L$ will be used in computed examples of Section \ref{sec4b}. 

We describe an approaches to compute a solution of \eqref{tikhonov2} and start with the 
simplest one. The matrix $M_k$ has a null space of large dimension (at least $n-k$). 
Therefore the Tikhonov minimization problem \eqref{tikhonov2} is not guaranteed to have a 
unique solution. To remedy this difficulty, we require the solution of \eqref{tikhonov2} 
to live in a subspace of fairly low dimension. A simple solution method is obtained when
using the solution subspace ${\mathcal R}(Q_k^{(r)})$. Then we obtain the minimization 
problem
\begin{equation}\label{tikhonov3}
\min_{\sby\in\R^k}\{\|M_k Q_k^{(r)}\by-\bg\|^2+\mu\|LQ_k^{(r)}\by\|^2\},
\end{equation}
which has a unique solution if and only if
\begin{equation}\label{nullcond2}
{\mathcal N}(M_k Q_k^{(r)})\cap{\mathcal N}(LQ_k^{(r)})=\{\bzero\}.
\end{equation}
This holds, in particular, when the triangular matrices $R_k^{(c)}$ and $R_k^{(r)}$ in 
\eqref{QRMk} are nonsingular. We found \eqref{nullcond2} to hold in all computed examples
that we solved.

Introduce the QR factorization
\begin{equation}\label{LQ}
LQ_k^{(r)}=Q_k^{(L)} R_k^{(L)},
\end{equation}
where the matrix $Q_k^{(L)}\in\R^{n\times k}$ has orthonormal columns and
$R_k^{(L)}\in\R^{k\times k}$ is upper triangular. We note that since the matrix $L$
typically is very sparse and $k$ is not large, the left-hand side of \eqref{LQ} generally
can be evaluated quite quickly also when $n$ is large. The minimization problem
\eqref{tikhonov3} yields the small problem
\begin{equation}\label{tikhonov4}
\min_{\sby\in\R^k}\{\|R_k^{(c)}(R_k^{(r)})^T\by-(Q_k^{(c)})^T\bg\|^2+
\mu\|R_k^{(L)}\by\|^2\}.
\end{equation}
This problem can be solved in several ways: We may compute a generalized SVD (GSVD) of the 
matrix pair $\{R_k^{(c)}(R_k^{(r)})^T,R_k^{(L)}\}$ (see, e.g., \cite{DNR,Ha1}), or apply a
cheaper reduction of the matrix pair that can be used when the generalized singular 
values of the matrix pair are not explicitly required; see \cite{DR}. 

When the matrix $R_k^{(L)}$ is nonsingular and not very ill-conditioned, which is the case
in many applications, one may consider transforming the minimization problem 
\eqref{tikhonov4} to standard form by the substitution $\bz=R_k^{(L)}\by$. This yields the
problem
\begin{equation}\label{reduced}
\min_{\sbz\in\R^k}\{\|R_k^{(c)}(R_k^{(r)})^T(R_k^{(L)})^{-1}\bz-(Q_k^{(c)})^T\bg\|^2+
\mu\|\bz\|^2\},
\end{equation}
which easily can be solved, e.g., by computing the singular value decomposition of the 
matrix $R_k^{(c)}(R_k^{(r)})^T(R_k^{(L)})^{-1}$. The solution $\bz_\mu$ of \eqref{reduced}
yields the solution $\by_\mu=(R_k^{(L)})^{-1}\bz_\mu$ of \eqref{tikhonov3}, from which we 
determine the approximate solution $\widetilde{\bx}_\mu=Q_k^{(r)}\by_\mu$ of 
\eqref{tikhonov}. The solution of \eqref{reduced} is cheaper than the solution of 
\eqref{tikhonov4} with the aid of the GSVD; see \cite{GVL} for counts of the arithmetic
floating point operations necessary to compute the GSVD of a pair of $k\times k$ matrices,
and the SVD of a $k\times k$ matrix.


\subsection{The discrepancy principle}\label{sec3_3}
We turn to the computation of the regularization parameter $\mu>0$ by the discrepancy 
principle. Assume for the moment that the matrix $A$ is available. Then we can solve 
equation \eqref{discrp} for $\mu>0$ by using a zero-finder such as Newton's method or one
of the zero-finders described in \cite{BPR,RS}. The theoretical justification of the 
discrepancy principle requires that the unavailable error-free vector $\widehat{\bg}$ 
associated with the available error-contaminated vector $\bg$ satisfies 
$\widehat{\bg}\in{\mathcal R}(A)$. 

Now consider the application of the discrepancy principle to the determination of the 
regularization parameter in \eqref{tikhonov2}. Generally, 
$\widehat{\bg}\not\in{\mathcal R}(M_k)$ and, therefore, the discrepancy principle cannot 
be applied when solving \eqref{tikhonov2} without modification. In the computed examples, 
we determine $\mu$ so that the computed solution ${\bx}_\mu$ of \eqref{tikhonov2} 
satisfies
\begin{equation}\label{discr2}
\|M_k\bx_\mu-Q_k^{(c)}(Q_k^{(c)})^T\bg\|=\eta\delta;
\end{equation}
cf. \eqref{discrp}. If the matrices $R_k^{(c)}$ and $R_k^{(r)}$ in \eqref{QRMk} are 
nonsingular, which generally is the case, then $Q_k^{(c)}(Q_k^{(c)})^T$ is an 
orthogonal projector onto ${\mathcal R}(M_k)$, and 
$Q_k^{(c)}(Q_k^{(c)})^T\widehat{\bg}$ lives in ${\mathcal R}(M_k)$. Equation 
\eqref{discr2} can be solved for $\mu\geq 0$ by using a zero-finder. 

\subsection{Stopping criterion for the adaptive cross approximation algorithm}
\label{sec3_4}
In view of \eqref{discrp}, we would like to determine a value of the regularization
parameter $\mu>0$ such that
\begin{equation}\label{discrp2}
\| A\bx_{\mu}(M_k) -\bg \| = \eta \delta.
\end{equation}
Even though the matrix $A$ is not available, we can determine an approximate upper bound 
for the left-hand side of \eqref{discrp2} as follows:
\begin{eqnarray}	
\nonumber
\|A\bx_{\mu}(M_k) - \bg \| &=& \|(A-M_k)\bx_{\mu}(M_k) + M_k\bx_{\mu}(M_k) - \bg \|\\
 \nonumber 
 &\leq& \|A - M_k\|\|\bx_{\mu}(M_k)\|+\| M_k\bx_{\mu} (M_k) - \bg \| \\
\label{eq900}
 &\lessapprox& S_k \| \bx_{\mu}(M_k)\| + \| M_k \bx_{\mu} (M_k) - \bg \|,
\end{eqnarray}
where $S_k$ is an approximation of $\| A - M_k\|$.  Based on the values of 
$(R_{k})_{i_\ell,j_\ell}$ from (\ref{Rk}) such an approximation $S_k$ can be computed as
\begin{eqnarray}
\nonumber
\| A - M_k\| &\leq& \|A-M_k\|_F \\
\nonumber
&\approx& \sqrt{\sum_{(i_\ell,j_\ell) \in \Pi} |(A-M_k)_{i_\ell,j_\ell} |^2} (mn) / |\Pi|\\
\nonumber
&=& \sqrt{\sum_{(i_\ell,j_\ell) \in \Pi}  |(R_k)_{i_\ell,j_\ell} |^2} (mn) / |\Pi| \\
\label{Sk}
&=:& S_k,
\end{eqnarray}
where $\|\cdot\|_F$ denotes the matrix Frobenius norm, and the set $\Pi$ is defined by
\eqref{set_pi} with $|\Pi| = t$ elements.

The number of step $k=k^*$ of the adaptive cross approximation algorithm is chosen to
be as small as possible such that there exists a $\mu>0$ such that
$$
S_k \|\bx_{\mu}(M_k)\| = \eta_1 \delta \mbox{ and } 
\| M_k \bx_{\mu}(M_k) - \bg \| \leq \eta_2 \delta.
$$
In this case we have, based on (\ref{eq900})
\begin{eqnarray*}
\|A \bx_{\mu}(M_k) - \bg \| &\lessapprox& S_k \| \bx_{\mu}(M_k)\| + 
\| M_k\bx_{\mu} (M_k) - \bg \| \\
&\lessapprox& (\eta_1 + \eta_2) \delta.
\end{eqnarray*}

\section{Numerical experiments}\label{sec4b}
This section describes a few computed examples with the adaptive cross approximation
method. For problems in one space-dimension, we will use the regularization matrices
$L_0 = I$, where $I\in\R^{n\times n}$ denotes the identity matrix, and
\begin{equation}\label{L1}
L_1 = \frac{1}{2}\left[ \begin{array} {cccccc}
 1 &   -1    &     &    &   & \mbox{\Large 0} \\
   &  \phantom{-}1    &   -1    &     &    & \\
   &   &  \phantom{-}1    &  -1  &    &    \\
   &   &  & \ddots & \ddots  &   \\
   \mbox{\Large 0}   &        &        &    & \phantom{-}1  & -1
   \end{array}
   \right]\in{\R}^{(n-1)\times n} 
\end{equation}
or
\begin{equation}\label{L2}
L_2 = \frac{1}{4}\left[ \begin{array} {cccccc}
-1   &  \phantom{-}2    &   -1    &     &    &\mbox{\Large 0} \\
&   -1  &  \phantom{-}2    &  -1  &    &    \\
&       & \ddots & \ddots & \ddots  &   \\
\mbox{\Large 0}    &        &        &   -1  & \phantom{-}2      & -1
\end{array}
\right]\in{\R}^{(n-2)\times n},
\end{equation}
which approximate a multiple of the first and second order derivative operators, 
respectively, assuming that $\bx_\mu$ is the discretization of a function 
$x_\mu$ at equidistant points on a bounded interval $\Omega_1$. For problems in two 
space-dimensions, we use the regularization matrices
\begin{equation}\label{Ltensor1}
L_{1,\otimes}=\left[\begin{array}{ccc} I & \otimes & L_1 \\ L_1 & \otimes & I 
\end{array}\right]
\end{equation}
or
\begin{equation}\label{Ltensor2}
L_{2,\otimes}=\left[\begin{array}{ccc} I & \otimes & L_2 \\ L_2 & \otimes & I 
\end{array}\right],
\end{equation}
where $\otimes$ stands for the Kronecker product. These choices of regularization matrices
are fairly common; see, e.g., \cite{Ha1,HRY,PRRY} for illustrations.

The examples are taken from the Regularization Tools \cite{PCH4} and from IR Tools 
\cite{GHN}. In all examples, the number of elements in the set $\Pi$ is $50 n$. The 
values of the parameters $\eta,\eta_1,\eta_2$ are chosen as $\eta = \eta_1=\eta_2 = 1.0$.

\begin{figure}
\centering
\includegraphics[scale=.45]{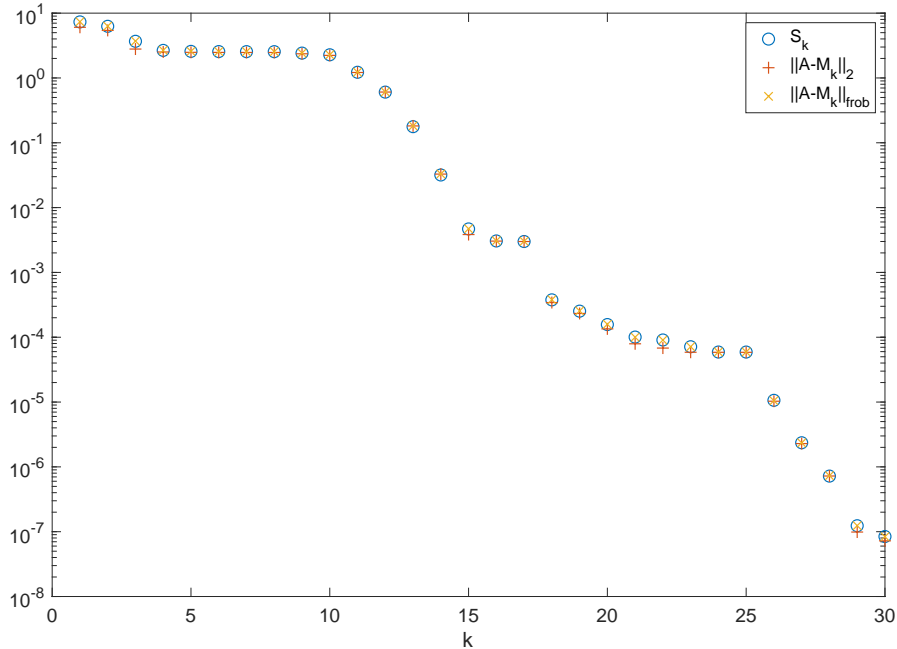}
\caption{For the ``gravity'' problem the values of $S_k$, $\|A-M_k\|$, and $\|A-M_k\|_F$
are plotted as a function of the number of steps $k$ of the adaptive cross approximation
algorithm.}\label{fig02}
\end{figure}

{\bf Experiment 1:}
We consider the problem ``gravity'' of size $n=1024$ from \cite{PCH4}. The aim of this 
experiment is to illustrate that the quantities $S_k$ defined by \eqref{Sk}  provide quite
accurate approximations of $\|A-M_k\|$. This is displayed by Figure \ref{fig02}.

\begin{figure}
\centering
\includegraphics[scale=.45]{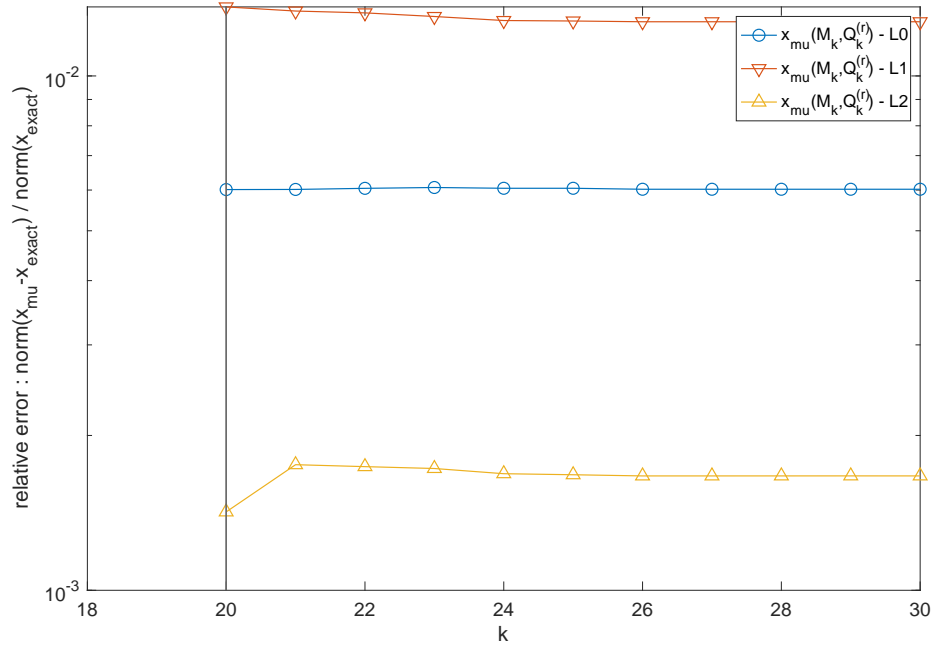}
\caption{For the ``gravity'' problem the relative error between the exact solution and the
computed approximate solution is plotted as a function of the iteration step $k$ and the 
chosen regularization matrix $L_0,L_1,L_2$. The size of the problem is $n = 1024$ with 
$\delta = 1.0e-2$ and $\eta = \eta_1= \eta_2=1.0$. The vertical line indicates where the 
stopping criterion is satisfied for adaptive rank approximation.}\label{fig01}
\end{figure}

{\bf Experiment 2:} 
We again consider the example ``gravity'' from \cite{PCH4} of size $n=1024$. Let
$\delta=10^{-2}$. This example illustrates that for certain problems only a fairly small 
number of steps of the adaptive cross approximation algorithm suffices to yield a 
satisfactory result. The example also shows that it may be very beneficial to use a 
regularization matrix different from the identity matrix. 
The maximum number of adaptive cross approximation steps is $30$. Results are shown for
$L\in\{L_0,L_1,L_2\}$. The quality of the computed solution $\bx_{\rm computed}$ is measured
by the relative error $\|\bx_{\rm computed}-\hat{\bx}\|/\|\hat{\bx}\|$. The horizontal
axis of Figure~\ref{fig01} shows the number of steps of adaptive cross approximation; the
vertical line indicates that for each one of the choices of $L$, the stopping criterion 
for the method is satisfied at step $k^* = 20$.

\begin{figure}
\centering
\includegraphics[scale=.45]{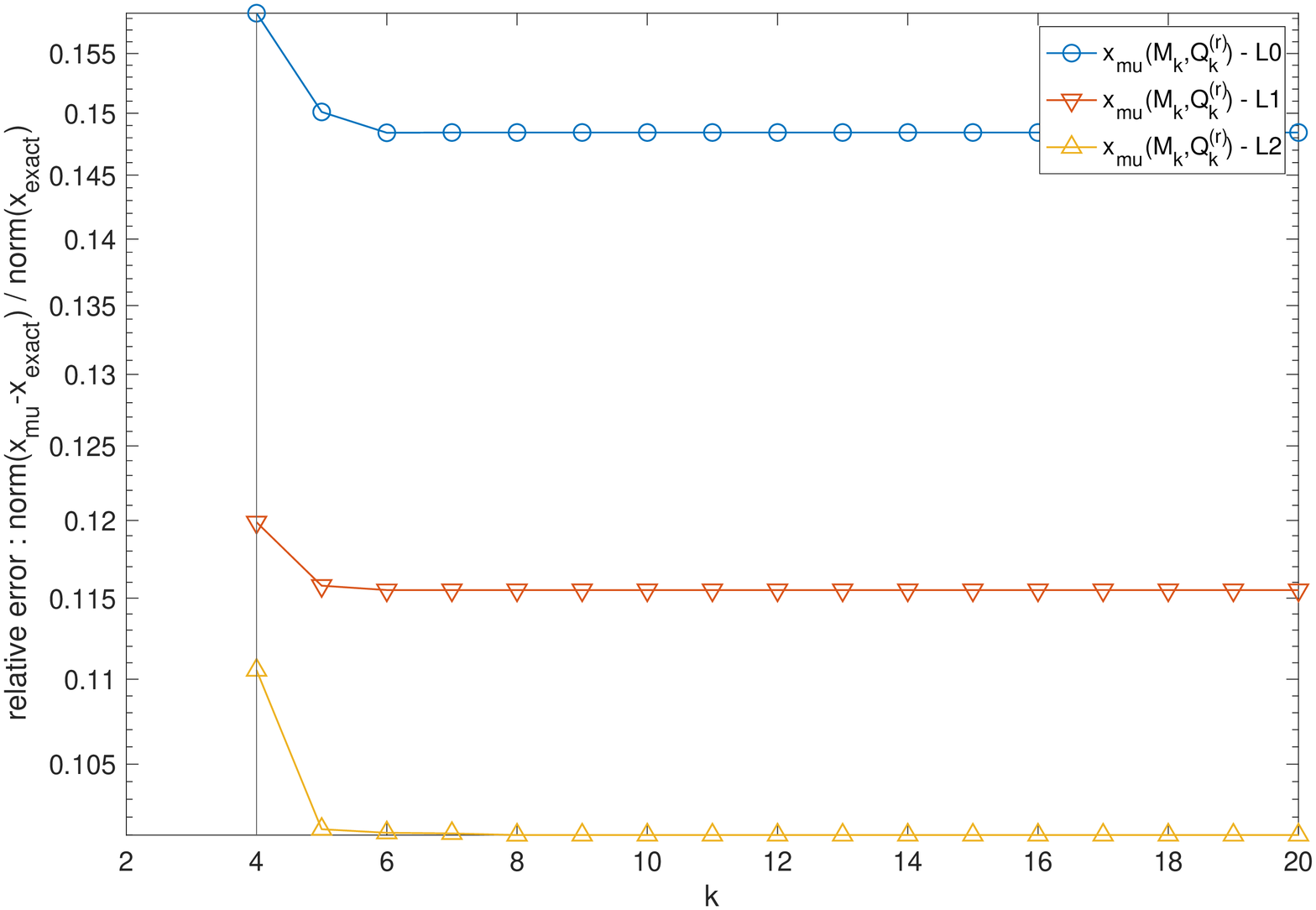}
\caption{For the ``baart'' problem the relative error between the exact solution and the 
computed approximate solution is plotted as a function of the iteration step $k$ and the 
chosen regularization matrix $L_0,L_1,L_2$. The size of the problem is $n = 1024$ with 
$\delta = 1.0e-2$ and $\eta =\eta_1=\eta_2= 1.0$. The vertical line indicates where the 
stopping criterion is satisfied for adaptive rank approximation.}\label{fig_baart01}
\end{figure}

{\bf Experiment 3:} This experiment is similar to Experiment 2, but for problem ``baart''
from \cite{PCH4}. Results are displayed in Figure~\ref{fig_baart01}. Also in this 
example it is beneficial to use a regularization matrix different from the identity.

\begin{figure}
\centering
\includegraphics[scale=.45]{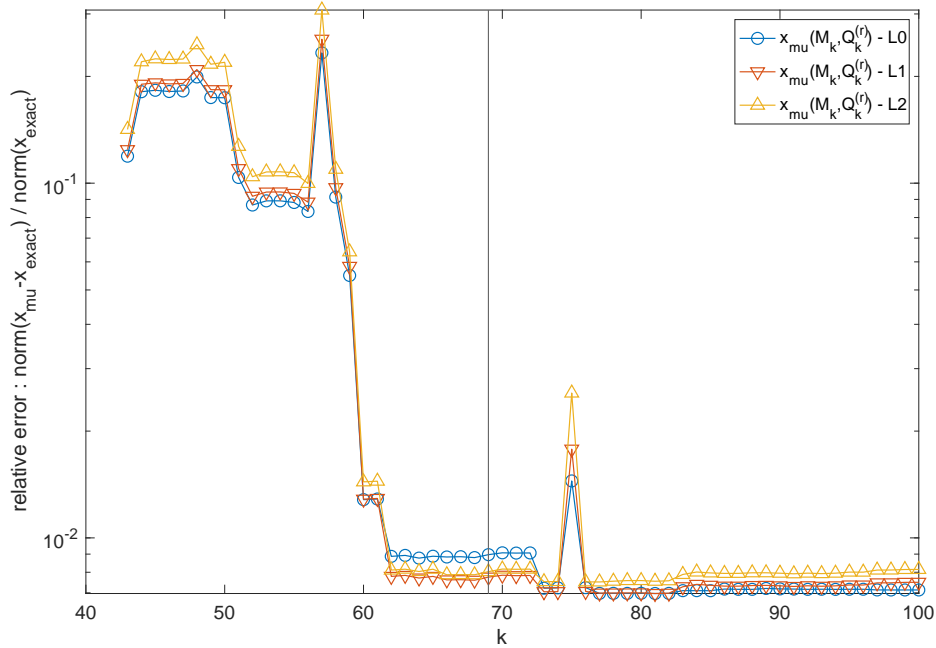}
\caption{For the ``phillips'' problem the relative error between the exact solution and 
the computed approximate solution is plotted as a function of the iteration step $k$ and 
the chosen regularization matrix $L_0,L_1,L_2$.  The size of the problem is $n = 1024$ 
with $\delta = 1.0e-2$ and $\eta = \eta_1=\eta_2=1.0$. The vertical line indicates where 
the stopping criterion is satisfied for adaptive rank approximation.}
\label{fig_phillips01}
\end{figure}

{\bf Experiment 4:} This expperiment is similar to Experiment 2, but for problem 
``phillips'' from \cite{PCH4}. The result is shown in Figure~\ref{fig_phillips01}.
For this example all three regularization matrices used perform about equally well.
The singular values of the matrix $A$ decay to zero slower for this example than for
the previous examples. Therefore more steps with the adaptive cross approximation 
algorithm have to be carried out.

%

{\bf Experiment 5:} We consider the example EXdiffusion\_rrgmres from the IR Toolbox
\cite{GHN}. The size of the problem is $4096$ and $\delta = 5.0e-3\|\widehat{\bg}\|$. The
other parameters are as in Experiment 2. The ``best'' solution determined by the example 
script in the IR Toolbox has relative error (when compared with the exact solution) 
$0.1875$, while for $L=L_2$ our algorithm in iteration step $k = 76$ reaches a relative 
error of $0.1910$.

\begin{figure}
\centering
\includegraphics[scale=.45]{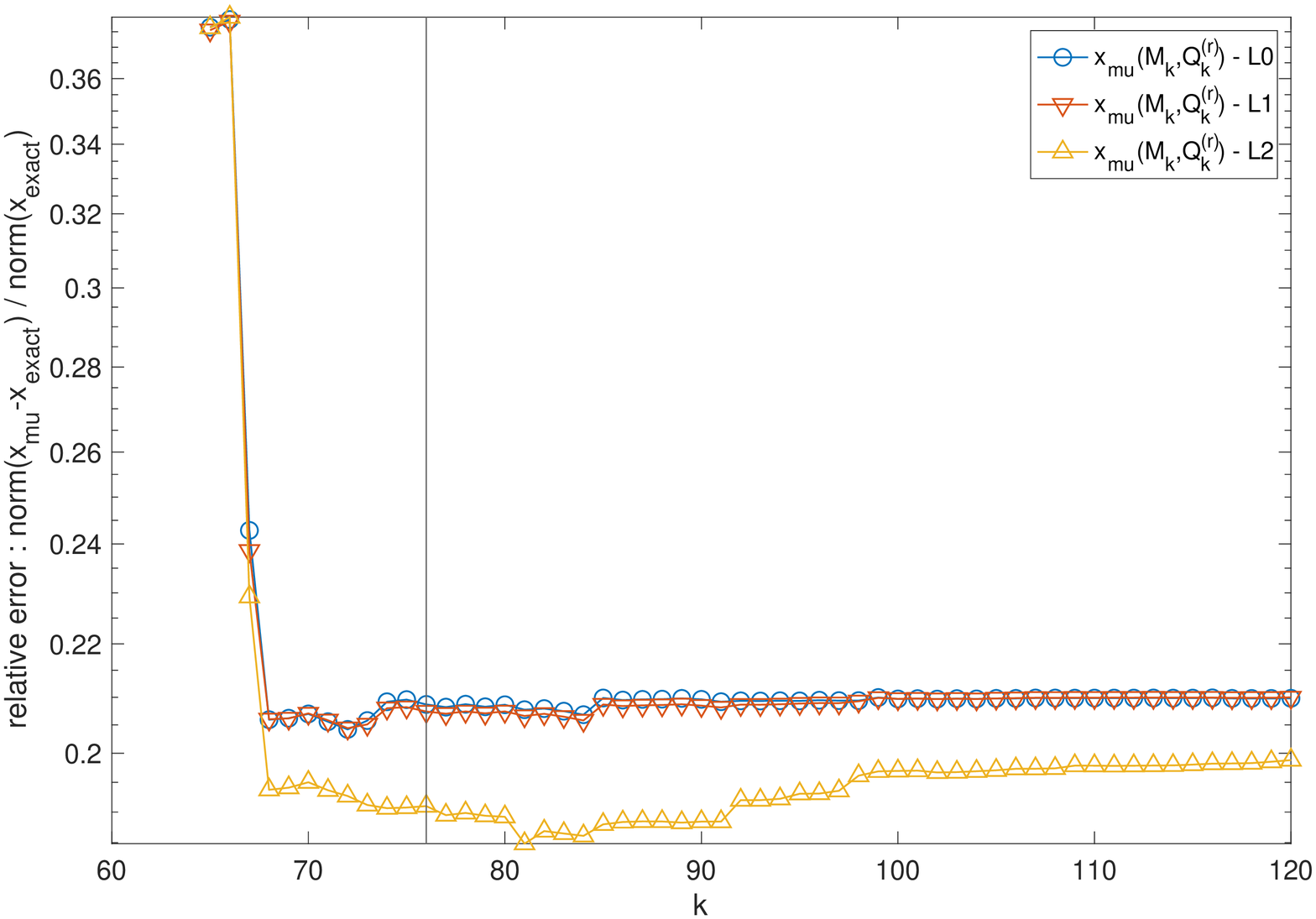}
\caption{For the ``diffusion'' problem the relative error between the exact solution and 
the computed approximate solution is plotted as a function of the iteration step $k$ and 
the chosen regularization matrix $L$. The size of the problem is $n = 4096$ with 
$\delta = 5.0e-3\| \widehat{\bg} \|$ and $\eta =\eta_1=\eta_2= 1.0$. The vertical line 
indicates where the stopping criterion is satisfied for adaptive rank approximation.}
\label{fig_diffusion}
\end{figure}

\begin{figure}
\centering
\includegraphics[scale=.45]{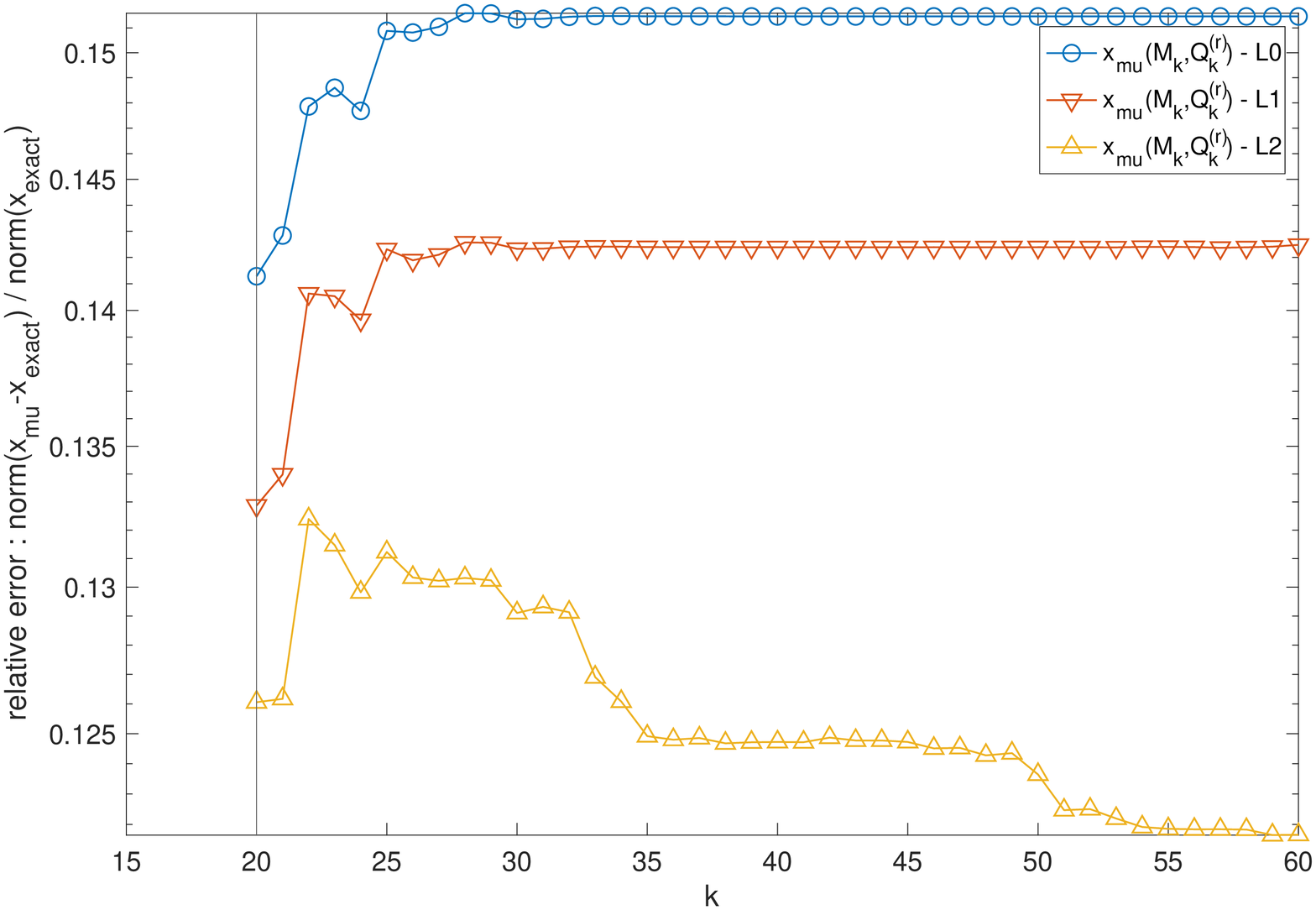}
\caption{For the 2D ``baart'' problem the relative error between the exact solution and 
the computed approximate solution is plotted as a function of the iteration step $k$ and 
the chosen 2D regularization matrix based on $L_0,L_1,L_2$. The size of the problem is 
$n = 40^2$ with $\delta = 1.0e-3$ and $\eta = \eta_1=\eta_2=1.0$. The vertical line 
indicates where the stopping criterion is satisfied for adaptive rank approximation.}
\label{fig_baart2D_02}
\end{figure}

{\bf Experiment 6:} The same as in Experiment 2 where we consider $A = B \otimes B$ with 
$B$ the matrix of the ``baart'' regularization problem. The matrix $B$ is of order $40$.
We use (\ref{Ltensor2}) for $L_0$, $L_1$ and $L_2$ as regularization matrices. The 
relative errors of the computed solutions are displayed in Figure~\ref{fig_baart2D_02}.

\section{Conclusion}\label{sec5}
This paper discusses the application of adaptive cross approximation to Tikhonov
regularization problems in general form. The computed examples illustrate that often 
only quite few cross approximation steps are required to yield useful approximate 
solutions. Particular attention is given to the stopping criterion for adaptive 
cross approximation.

\section*{Acknowledgment}
This research was partially supported by 
 the Fund for Scientific Research–Flanders (Belgium),  Structured Low-Rank Matrix/Tensor Approximation: Numerical Optimization-Based Algorithms and Applications: SeLMA, EOS 30468160,  
 the KU Leuven Research Fund,  Numerical Linear Algebra and Polynomial Computations, OT C14/17/073.

\end{document}